 \documentstyle[12pt]{article}
\newcommand{\beq}{\begin{equation}}
\newcommand{\eeq}{\end{equation}}
\newcommand{\beeq}{\begin{eqnarray}}
\newcommand{\ba}{\begin{array}{c}}
\newcommand{\ea}{\end{array}}
\newcommand{\eeeq}{\end{eqnarray}}
\newcommand{\bdm}{\begin{displaymath}}
\newcommand{\edm}{\end{displaymath}}

\newcommand{\ra}{\rightarrow}

\def\cs{coherent states }
\def\GL2{$GL_q(2)$ }
\def\ha{Hopf algebra }

\def\o{oscillator }

\def\q{quantum }
\def\Q{Quantum }
\def\qa{quantum algebra }

\def\YBe{Yang--Baxter equation }

\font\BBm=msbm10   scaled\magstep1
\font\got=eufm10 scaled\magstep1
\font\xx=eurb10

\begin{document}


\vspace{2.5cm}

\begin{center}

{\large {\bf Some Remarks on the Representations}}\\
{\large {\bf of the Generalized Deformed}} \\
{\large {\bf Oscillator Algebra}} \\[1cm]

{\bf V.V.Borzov } \\
 St Petersburg University of Telecommunications, \\
Department of Mathematics \\
Moika 61, 191065, St Petersburg, Russia \\[.5cm]

{\bf E.V.Damaskinsky} \footnote
 {  Supported by Russian Foundation for
 Fundamental Research, Grant N 95-01-00569-a}  \\
 Defense Constructing Engineering Institute, \\
 Zacharievskaya st 22, 191194, St Petersburg, Russia \\[.5cm]

{\bf S.B.Yegorov}  \footnote{Supported by Russian Foundation for
 Fundamental Research, Grant N 95-01-00569-a}
\footnote{ E-mail address: ye@stu.neva.ru  } \\
St.Petersburg Technical University. Computer Center.\\
Politechnicheskaya st 29, 195251,  St Petersburg, Russia \\ [1.2cm]

  {\bf Abstract.}
\end{center}

The classification of the representations of the generalized deformed
oscillator algebra is given together with several comments about possibility
of introducing  a coproduct structure in some type of deformed
oscillator algebra.

\newpage

\section{Introduction}
In recent years a lot of interest has been devoted to the study of the
various quantum deformations of boson oscillator commutation relations
algebra (see f.ex.\cite{4}). From mathematical point of view such popularity
connected with numerous relations which exsist between deformed
oscillators and other quantum deformations (quantum groups, quantum algebras,
quantum spaces etc \cite{*3}-\cite{*6}) which was born in the course of
development of the quantum inverse scatering method \cite{*1}-\cite{*2}.
{}From the other side there are some hopes that in a physical studies of
non-linear phenomena the deformed oscillator can play the role much the
same as the usual boson oscillator in standard quantum mechanics. Such
hopes are supported by several applications of the deformed oscillators
in conformal field theory \cite{*7}, lattice models \cite{*8,*9}, nuclear
spectroscopy \cite{*10,*11}, in describing the systems with non-standard
statistics and energy spectrum \cite{BEN} etc.

\vspace{12pt}

In this work we describe all irreducible representations of the so-called
{\bf generalized deformed oscillator algebra} (GDO-algebra) \cite{66} which
contains as particular cases the most popular variants of such deformations
together with the standard (undeformed) one. We also give some comments
connected with possibility of introducing the Hopf algebra structure.
In particular we show that to this end we must essentially
restrict the algebra for example by introducing in the definition of the
algebra a mate partner to the main commutation relation .  Such extension of
defining relations has grave consequences.  Namely, modificated  algebra
admit much less representations: as a rule only representation similar to
Fock one are survived. Such "uniqueness" is probably not very desirable
because there are some indications (see f.ex. \cite{BEN}) of possible
physical significance of non Fock representations.

In the rest of this section mainly for purpose of establishment notations for
future reference we recall briefly some well known facts about standard
boson oscillator algebra essential for the following discussions.

\vspace{8pt}

Let us recall that the {\bf Heisenberg algebra ${\cal H}$} is a Lie
algebra, generated by creation operator $a^+$ and annihilation operator $a$
and unity, subject to the commutation rules
\beq\label{a1}[a,a^{+}]\equiv aa^+-a^+a=I,\quad [a,I]=0,\quad [a^+,I]=0.\eeq
The related Lie group is simplest nilpotent group, known as Heisenberg group
(see for example \cite{*12}-\cite{*14}). The number
operator is defined in this case by the relation \beq\label{a2} N:=a^+a.\eeq
This operator fulfil the commutation rules (the Heisenberg "equations of
motions") \beq\label{a3}[N,a]=-a,\qquad [N,a^{+}]=a^{+},\eeq with generators
$a^+$ and $a$. Note that the number operator $N$ is not included in the list
of generators for ${\cal H}$. This is essential remark for
understanding some of differences in definitions of various deformations of
relation (\ref{a1}).

It is well-known that in accordance with von Neumann uniqueness
theorem (see for example (\cite{*11,*14})) the Heisenberg group has
essentially only one (up to unitary equivalence) irreducible representation.
The similar, but much weaker uniqueness assertion \cite{*14} is true for
representations of the Heisenberg algebra.

In the form of the {\bf Fock} (or {\bf occupation number})
{\bf representation} it can be described in the Dirac bracket notation
by the following formulas. We denote by $|0>$ the (unique) vacuum state which
is the solution of the equation $a|0>=0$ and construct the other basis states
in standard way
\beq \label{b1}|n>={\frac {(a^{+})^n}{\sqrt{n!}}}|0>,\quad n\in
{\hbox{\BBm N}}.\eeq The Fock Hilbert space {\xx H}$_F$ is defined by
properly completing of the linear span of these vectors in the scalar product
defined as $<m|n>=\delta _{m,n}$. Then the Fock representation is defined by
the following action  on this basis states of the main operators
\beq \label{b2}\begin{array}{l}
a|n>=\sqrt{n}|n-1>, \; n\geq 1;\quad a|0>=0; \\
a^+|n>=\sqrt{n+1}|n+1>;\qquad N|n>=n|n>. \ea \eeq

The quite similar but different algebraic object is the {\bf boson
oscillator algebra ${\cal W}^0_{0,0}(1)$}.  This is a Lie algebra generated
by {\bf three} generators $\{a,a^+,N\}$ and unity, which have the following
commutation rules \beq \label{a4} \ba \lbrack a,a^{+}\rbrack=I,\qquad
\lbrack I,all \rbrack=0 \\ \lbrack N,a \rbrack =-a,\qquad  \lbrack N,a^{+}
\rbrack =a^{+}.\ea\eeq The related Lie group is solvable group, known as
oscillator group (see for example \cite{d3}).

We would like to stress that in the case of the oscillator algebra
${\cal W}^0_{0,0}(1)$ the operator $N$ is an independent generator and
commutation rules (\ref{a3}) are now postulated and they are not  deduced as
above. This means that the connection (\ref{a2}) is not supposed and indeed
it does not hold for all representations. In contrast to the Heisenberg
case, the oscillator group (and the oscillator algebra too) has reach enough
supply of nonequivalent irreducible representations. Therefore the relation
(\ref{a2}) is considered in this case as constraint for selection of
representations with needed properties.

Let us note, that the question about compatibility between commutation
relations (\ref{a1}) and "equations of motion"  (\ref{a3}) was raised by
Wigner in \cite{d5}, and studied in some details in \cite{d6}-\cite{d7}

As pointed in \cite{4,63} the mentioned distinction between these two
algebras are important to have in a mind when we consider the different
deformations of the boson oscillator especially in discussion of
representations and possible Hopf algebra structure.

\section{ Generalized Deformed Oscillator Algebra
${\cal W}^{\gamma}_{\alpha,\beta}(q)$}
In this section we list some of the most popular types of the  possible
deformations and unify them in the unique generalized deformed oscillator
algebra which was slightly generalized variant of the algebra suggested in
\cite{66}.

All these deformations have some common peculiarities. First, such
algebras are not Lie algebras and rigorously defined as quotient
${\cal A}_q={\cal A}/I_q$ of associative algebra ${\cal A}$ by two-sided
ideal $I_q$. Here the associative algebra ${\cal A}$, freely generated by
fixed set of letters (generators), and ideal $I_q$ defined by fixed set of
relations (commutation rules).

Second, in such algebras any relations of the type (\ref{a2}) are in
general not hold and may appear only as constraints on representations, or
must be postulated independently in the definition of the algebra.

Third, it is always assumed the some kind of correspondence rule. Name\-ly,
it was assumed that in the so-called "classical limit" $q\ra 1$ the
structure of Lie algebra of one of the usual types mentioned in Introduction
is retained.  But nevertheless this it is possible appearing of the so called
singular representations \cite{6} which are absent in the such limit is
possible .

{\bf Tamm-Dancoff oscillator algebra ${\cal W}^{1}_{1,0}(q)$}.
This algebra arised in the frames of the Tamm-Dancoff method \cite{68},
\cite{69} in quantum field theory and based on the commutator rules
\beq\label{a5}aa^+-qa^+a=q^N, \quad  [N,a]=-a,\quad [N,a^{+}]=a^{+}, \eeq
where $q$ is an arbitrary (in general complex) non-zero number.

{\bf Arik-Coon-Kuryskin oscillator algebra
${\cal W}^{1}_{0,0}(q)$}. This algebra was introduced independently by many
authors (see for example \cite{8}-\cite{40}). Defining commutation rules for
it looks as \beq\label{a6}aa^+-qa^+a=I,\eeq which can be extended by
(\ref{a3}) as additional relations. Of course in this case the equation
(\ref{a2}) is not hold.

{}From the relation (\ref{a6}) we can easely receive the useful relation
\beq\label{a7}a(a^+)^m-(qa^+)^ma=[m;q](a^+)^{m-1},\quad [m;q]\equiv
      {\frac {q^m-1}{q-1}}.   \eeq
We note that appearance of the basic number $[m;q]$ in this formula
indicates the connection of this algebra with well developed mathematical
theory known as basic analysis (or q-analysis) and allows to apply in the
studying of its representations the powerful methods of the theory of basic
hypergeometric series \cite{*15}.

This algebra has nontrivial center \cite{4} generated by
\beq\label{a8} \zeta= [N;q]-a^+a.  \eeq
Nontriviality of the center explained the existence of non equivalent
irreducible representations  of this algebra (absence of the counterpart
of unique\-ness theorem for Heisenberg algebra) \cite{6},\cite{4},
\cite{32}-\cite{35}.

Attractive peculiarities of this algebra consist in that it is possible to
define in more or less simple way the multi-mode generalization of it in
non trivial (not mutually commuting) fashion which is covariant under the
quantum group $SU_q(2)$ \cite{32}-\cite{33},\cite{71}.

{\bf Quantum deformed oscillator algebra ${\cal W}^{1}_{-1,0}(q)$}.
This algebra defined by the relations
\beq\label{a9}aa^+-qa^+a=q^{-N},\quad [N,a]=-a,\quad [N,a^{+}]=a^{+}.  \eeq
It aroused in attempts to extend the wellknown Schwinger boson realisation
\cite{d8} of angular momentum operators for the case of quantum algebras
$U_q(su(2))$ \cite{3,2} and its non-compact form $U_q(su(1,1))$ \cite{44,4}.

{}From (\ref{a9}) we can receive the following useful formulas \cite{4}
\beq\label{a10} a^+q^{-N}=qq^{-N}a^+,\qquad aq^{-N}=q^{-1}q^{-N}a; \eeq
\beq\label{a11} [N,a^+a]=0,\qquad [N,aa^+]=0; \eeq
\beq\label{a12} a(a^+)^m-(qa^+)^ma=[m](a^+)^{m-1}q^{-N},   \eeq
where
\beq\label{a14} [m]\equiv {\frac {q^m-q^{-m}}{q-q^{-1}}}=
  {\frac {{\hbox{\rm sh}}(\eta m)}{{\hbox{\rm sh}}(\eta)}},\quad
  {\hbox{\rm where}}  \;q=e^{\eta}. \eeq

Thus we see that this algebra connected with the symmetrical under the
replacement $q\ra q^{-1}$ basic number $[m]$. Unfortunately we have not the
developed variant of $q$-analysis based on this kind of the basic number.

This algebra also has nontrivial center \cite{4} generated by the element
similar to (\ref{a8}) and thus there are many non-equivalent irreducible
representations also for this algebra \cite{7},\cite{6},\cite{4}.

Some authors considered restricted form ${\cal W}^{1}_{-1,0}(q;q^{-1})$
of this algebra in which besides the relations (\ref{a9}) the additional
relation  \beq\label{a15} aa^+-{\frac 1q}a^+a=q^{N},  \eeq
is  postulated too. In this case we also have
\beq\label{a16} a^+a=[N],\qquad aa^+=[N+1],\qquad [a,a^+]=[N+1]-[N]  \eeq
and as follows from (\ref{a8}) (with $[N;q]$ replaced by $[N]$) the center
of the algebra became trivial ($\zeta =0$). For this reason the resticted
algebra has only one (up to equivalence) irreducible representation, which
is the $q$-analogue of the standard Fock representation of the usual boson
oscillator.

Let us note that in the Fock representation we have direct connection
between usual non deformed operators $b, b^+$ and $N_b=b^{+}b$,
and deformed ones. This connection is given \cite{4} by
\beq N_F=N_b, \hspace{.2cm}
a_F^+ = \sqrt{ {[N_F]}\over{N_F} }b^+ , \hspace{.2cm}
a_F = \sqrt{{[N_F + 1]}\over{N_F}}b. \eeq

Provided that the this map connected the operators of deformed
$q$-oscilla\-tor with usual one is invertible (this is not the case when
$q^M=1$, \,$M \in ${\bf N}) this restricted algebra
${\cal W}^{1}_{-1,0}(q;q^{-1})$  is equivalent with standard quantum
mechanical boson oscillator algebra ${\cal W}^{0}_{0,0}(1)$. On the
other hand the algebra ${\cal W}^{1}_{-1,0}(q;q^{-1})$ can be identified
with $sl_q(2)$ {\qa}. Indeed if both relations
\beq aa^{+} -qa^{+}a=q^{-N}, \hspace{.5cm}
aa^{+} -q^{-1}a^{+}a=q^{N}, \eeq are valid, then operators
\beq X_+=\sigma a, \hspace{.3cm} X_-=\sigma a^\dagger ,
\hspace{.3cm} J=1/2(N-{{\pi i}\over{2\eta}}), \hspace{.5cm} (q=e^\eta) \eeq
where ${\sigma}^2={{i\sqrt{q}}\over{q-1}}$, fulfil commutation relations of
$sl_q(2)$. This equivalence of course allows one to induce the {\ha}
structure (similar to given in \cite{Y})  in ${\cal W}^{1}_{-1,0}(q;q^{-1})$
from $sl_q(2)$, but corresponding co-product do not respect Hermitian
conjugation.

{\bf Yan oscillator algebra ${\widehat W}(q)$ } This algebra was introduced
in the work \cite{Y} and based on the relations (\ref{a16}) which are not
equivalent with (\ref{a9}) and (\ref{a15}) (see discussion in \cite{63}).
As shown in \cite{Y} Hopf algebra structure for this algebra  can be
given. We note that this algebra has trivial center and only Fock-type
representations.

{\bf Feinsilver oscillator algebra ${\cal W}^{0}_{-2,0}(q)$}.  This
algebra defined in \cite{30}-\cite{31},\cite{4} by the relations
\beq\label{a17}[a,a^+]=q^{-2N},\quad [N,a]=-a,\quad [N,a^{+}]=a^{+}. \eeq

We note that algebra ${\cal W}^{0}_{-2,0}(q)$ can be obtained by contraction
to limit of infinite spin of the quantum algebra $U_q(sl_2)$
\cite{43}-\cite{36},\cite{4}. The center of ${\cal W}^{0}_{-2,0}(q)$ is
generated by
\beq\label{a18}\zeta=[N;q^{-2}]-a^+a,   \eeq
and some of its representations are considered in \cite{31}.

This algebra conclude the list of examples which we intend to describe here.
We note that more details about properties of these algebras and connections
exsisting between them may be founded for example in \cite{4},\cite{45},
\cite{36}.

Now we are ready to consider the main object of our discussions the algebra
which in particular cases gives each algebra described above (except the
${\widehat W}(q)$).
...

{\bf The generalized deformed oscillator algebra
${\cal W}^{\gamma}_{\alpha,\beta}(q)$}.
The generators of this algebra fulfil commutation relations
\beq\label{a19}aa^+-q^{\gamma}a^+a=q^{\alpha N+\beta},\quad
       [N,a]=-a,\quad [N,a^{+}]=a^{+},  \eeq
where $\alpha,\,\beta,\,\gamma$ are real parameters. From first of this
relations we receive by induction
\beq\label{a20}a(a^+)^m-q^{m\gamma}(a^+)^ma=(a^+)^{m-1}q^{\alpha
N+\beta}{\frac {q^{m\alpha}-q^{m\gamma}}{q^{\alpha}-q^{\gamma}}},\eeq
which is useful in computations.

Let us consider the (abstract) Fock representation of this algebra. To this
end we suppose that there exist the vacuum state $|0>$ such that
\beq\label{a21} a|0>=0,\qquad N|0>=0 \eeq
Then as in \cite{66} by standard technique
 we receive
\beq\label{a22}\begin{array}{l}  N|n>=n|n>,\\ \smallskip \\
    a|n>=\sqrt{F^{\gamma}_{\alpha,\beta}(n;q)}|n-1> ,\\ \smallskip \\
    a^+|n>=\sqrt{F^{\gamma}_{\alpha,\beta}(n+1;q)}|n+1>, \ea \eeq
on the states \beq\label{a23} |n>={\frac {(a^+)^n|0>}
{\sqrt{F^{\gamma}_{\alpha,\beta}(n;q)!}}} \quad n=0,1,2\dots, \eeq where
generalized basic number defined by \beq\label{a24}\ba
F^{\gamma}_{\alpha,\beta}(n;q)= \cases{ q^{\beta}{\frac
     {q^{n\alpha}-q^{n\gamma}}{q^{\alpha}-q^{\gamma}}} &
      $\alpha \neq \gamma$, \cr
     nq^{\beta +\gamma(n-1)} & $\alpha =\gamma $. \cr } \\ \smallskip \\
      F^{\gamma}_{\alpha,\beta}(n;q)!= F^{\gamma}_{\alpha,\beta}(1;q) \cdot
     F^{\gamma}_{\alpha,\beta}(2;q) \cdot \dots \cdot
      F^{\gamma}_{\alpha,\beta}(n;q),\quad
      F^{\gamma}_{\alpha,\beta}(1;q)=q^{\beta}.\ea	  \eeq

We note that for special values of parameters we received, respectively
\beq\label{a25}F^{1}_{-1,0}(n;q)={\frac {q^{n}-q^{-n}}{q-q^{-1}}} =
    [n]_q=[n],\eeq
which gives the symmetrical basic number (\ref{a14}) for quantum deformed
oscillator algebra ${\cal W}^{1}_{-1,0}(q)$. Similarly,
\beq\label{a26}F^{1}_{0,0}(n;q)={\frac {q^{n}-1}{q-1}} =[n;q], \eeq
that is standard basic number (\ref{a7}), related to Arik-Coon-Kuryskin-
Jannussis-Cigler oscillator algebra $\!{\cal W}^{1}_{0,0}(q)$.
For Tamm-Dancoff oscillator algebra $\!{\cal W}^{1}_{1,0}(q)$ we get
\beq\label{a27}F^{1}_{1,0}(n;q) =nq^{n-1}, \eeq
as the related basic number. Feinsilver oscillator algebra
${\cal W}^{0}_{-2,0}(q)$ related with
\beq\label{a28}F^{0}_{-2,0}(n;q)=q^{2(1-n)}[n;q^2]=
q^{2(1-n)}(1+q^n+q^{2n}).\eeq

Let us list some of the properties of this basic number which we used below:

\begin{equation}
\label{f7}\left\{
\begin{array}{ccc}
n<0 & \Rightarrow & \left\{
\begin{array}{ccc}
F_{\alpha ,\beta }^\gamma (n)<0 & if & \gamma \neq 0 \\
F_{\alpha ,\beta }^0(n)>0 & if & \gamma =0
\end{array}
\right. \\
n\geq 0 & \Rightarrow & F_{\alpha ,\beta }^\gamma (n)>0
\end{array}
\right. ;
\end{equation}
\begin{equation}
\label{f8}F_{\alpha ,\beta }^\gamma (n+1)-q^\gamma F_{\alpha ,\beta }^\gamma
(n)=q^{\alpha n+\beta };
\end{equation}

where $F_{\alpha ,\beta }^\alpha (n)=F_{\alpha ,\beta }^\alpha (n;q)$.
Finally, in "classical limit" we obtain naturally
\beq\label{a29}\lim\limits_{q\ra 1}F^{\gamma}_{\alpha,\beta}(n;q)=n. \eeq

\section{Classification of the irreducible representations of the
generalized \hfil \break deformed oscillator algebra
${\cal W}_{\alpha ,\beta}^\gamma (q) $}

In this section we consider generalized deformed oscillator algebra ${\cal
W}={\cal W}_{\alpha ,\beta }^\gamma (q)$ with defining relations
\begin{equation}
\label{f1}aa^{+}-q^\gamma a^{+}a=q^{\alpha N+\beta },
\end{equation}
\begin{equation}
\label{f2}[N,a]=-a,\qquad [N,a^{+}]=-a^{+}
\end{equation}
where parameters $\alpha ,\beta ,\gamma \in {\hbox{\BBm R}},$ and
$q$ is an arbitrary positive
number. Let operators $a,a^{+},N$ gives the realization of these relations in
separable Hilbert space ${\cal H}$ and fulfil natural hermiticity conditions
\begin{equation}
\label{f3}(a)^{*}=a^{+},\quad (a^{+})^{*}=a,\quad (N)^{*}=N.
\end{equation}
We also supposed that selfadjoint (Hermitian) operator $N$ has the simple
(non degenerate) discrete spectrum.

Here we give, following mainly the Rideau work
\cite{7} the classification of the irreducible representations of the
generalized deformed oscillator algebra ${\cal W}_{\alpha ,\beta }^\gamma
(q).$ Let us note that the selection rule for testing the type of the given
representation, suggested in \cite{34}, also
works in the more general case of the algebra ${\cal W}_{\alpha ,\beta
}^\gamma (q)$ considered here. But naturally we deal with modified selector
operator \begin{equation} \label{f4}K:=aa^{+}-q^\alpha a^{+}a.
\end{equation}

In the following we need to consider the five different domains for
parameters $\alpha,\;\gamma$ of the algebra
${\cal W}_{\alpha ,\beta }^\gamma (q):$
\begin{equation}\label{f5}\begin{array}{ccccc}
1)\,\gamma \geq 0,\;\alpha <\gamma \,; & \quad &
2)\,\gamma \geq 0,\;\alpha >\gamma \,; & \quad &
3)\,\gamma \leq 0,\;\alpha <\gamma \,; \\
4)\,\gamma \leq 0,\;\alpha >\gamma \,; & \quad &
5)                \;\alpha =\gamma \,; & \quad &
\end{array}\end{equation}

We note that parameter $\beta$ plays no role in the classification of the
representations for the algebra ${\cal W}_{\alpha,\beta }^\gamma (q):$.
Therefore the case 4) reduced to the case 1) (and the case 3) to the case 2)
too) with the help of the substitution $q\leftrightarrow q^{-1}.$

The central element of the algebra ${\cal W}_{\alpha ,\beta }^\gamma (q)$
can be written as
\begin{equation} \label{f10}\widehat{\zeta }=\left( F_{\alpha ,\beta}^\gamma
(N)-a^{+}a\right) {\cal S}(N),\end{equation}
where
\begin{equation}
\label{f11}{\cal S}(N)=\sum_{k=0}^\infty v_kN^k
\end{equation}
with arbitrary real coefficients $v_k.$

Indeed, because from (\ref{f2}) it follows that $[N,a^{+}a]=0,$ we have that
if $N\psi _0=\nu _0\psi _0,$ then $a^{+}a\psi _0$ is also the eigenvector
for $N$ with the same eigenvalue $\nu _0.$ Taking into account that $N$ has
the simple spectrum we obtain
\begin{equation}
\label{f12}a^{+}a\psi _0=\lambda _0\psi _0,
\end{equation}
where $\lambda _0\geq 0$ because of the operator $a^{+}a$ is non-negative.
It is standard task to check that the vectors ($a^{+})^n\psi _0$ and
$a^m\psi_0,$ if non zero, are the eigenvectors for $N$ with the eigenvalues
$\nu_0+n $ and $\nu _0-m,$ respectively$.$ Then these vectors
($a^{+})^n\psi _0$ and $a^m\psi _0,$ are also the eigenvectors for
$a^{+}a$ with the eigenvalues $\lambda _n$ and $\lambda _{-m}$, respectively.
Analogously
they are the eigenvectors for $aa^{+}$ with the eigenvalues $\mu _n$ and
$\mu _{-m}.$

Acting by the relation (\ref{f1}) on
\begin{equation}
\label{f13a}\psi _n=\left\{
\begin{array}{lll}
(a^{+})^n\psi _0 & for & n\geq 0 \\
a^{-n}\psi _0 & for & n<0
\end{array}
\right.
\end{equation}
we receive
\begin{equation}
\label{f13}\mu _n-q^\gamma \lambda _n=q^{\alpha (n+\nu _0)+\beta }.
\end{equation}
But for all $n\in {\hbox{\BBm Z}}$ such that $\psi _n\neq 0$ we have $\mu
_n=\lambda _{n+1}.$ So from (\ref{f13}) it follows the recurrent relation
\begin{equation}
\label{f15}\lambda _{n+1}=q^\gamma \lambda _n+q^{\alpha (n+\nu _0)+\beta },
\end{equation}
which has the solution
\begin{equation}
\label{f16}\lambda _n=q^{\gamma n} \lambda _0+q^{\alpha \nu _0}F_{\alpha ,\beta
}^\gamma (n).
\end{equation}
We need only those solutions of (\ref{f16}) for which $\lambda _n\geq 0$
(because of the operator $a^{+}a$ is non-negative). According to (\ref{f7})
the rhs (\ref{f16}) is always positive when $n\geq 0$ . But for $ n<0\;$
we have $\lambda _n\geq 0$ only if \begin{equation} \label{f17}\lambda
_0\geq -q^{\alpha \nu _0-\gamma n}F_{\alpha ,\beta }^\gamma (n),\quad n<0.
\end{equation} It is convenient to rewrite this inequality in the form
\begin{equation} \label{f17'}\lambda _0\geq \left\{ \begin{array}{lll}
q^{\alpha \nu _0+\beta -\alpha }\frac{1-q^{n(\alpha -\gamma )}}{1-q^{\gamma
-\alpha }}, & if & \alpha \neq \gamma  \\ -nq^{\gamma (\nu _0-1)+\beta }, &
if & \alpha =\gamma \end{array} \right. ;\quad (n<0) \end{equation}

Note that if
\begin{equation}\label{f18}q>1,\;\alpha \leq \gamma \quad {\rm or}\quad
0<q<1,\;\alpha \geq \gamma \end{equation}
then rhs of (\ref{f17}) goes to $+\infty $ as $n\rightarrow (-\infty ).$
Thus there exist $n_0<0$ such that $\lambda _n\leq 0,$ for all $n<n_0.$

On the other hand if \begin{equation}\label{f19}q>1,\;\alpha >\gamma \quad
{\rm or}\quad 0<q<1,\;\alpha <\gamma \end{equation}
the inequality (\ref{f17}) is fulfilled for all $n<0$ if
\begin{equation}
\label{f20}\lambda _0\geq \frac{q^{\alpha (\nu _0-1)+\beta }}{1-q^{\gamma
-\alpha }}.\end{equation}
In the opposite to (\ref{f20}) case there exist $n_0<0$ such that
$\lambda _n\leq 0,$ for all $n<n_0.$

Now we are ready to consider the irreducible representations of
the algebra ${\cal W}_{\alpha ,\beta }^\gamma (q)$ for the each case listed
in (\ref{f5}).

\subsubsection{1) The case $\gamma \geq 0,\quad \alpha <\gamma $.}

We suppose that the following inequalities
\begin{equation}\label{f20a}{\bf (A)\qquad }q>1\quad {\rm or}\quad 0<q<1,\;
\lambda _0<\frac{q^{\alpha (\nu _0-1)+\beta }}{1-q^{\gamma -\alpha }}.
\end{equation}   are hold.

According to (\ref{f16}) if $n \ra +\infty$ we have
\begin{equation}\label{f20a1}\begin{array}{ll}
\lambda _n \ra +\infty & {\rm if}\;\cases{
		    \gamma >0,\;     q>1,  & or \cr %
\noalign{\smallskip}
		    \gamma =0,\; 0< q <1,  & or \cr 
         \gamma >0,\;\alpha <0,\; 0< q <1;  &    \cr}
\\ \smallskip \\
\lambda _n \ra 0^+     & {\rm if}\; \gamma >0,\;\alpha >0,\; 0< q <1;
\\ \smallskip \\
\lambda _n \ra \lambda _0+\textstyle{\frac {q^{\alpha \nu _0 +\beta}}
       {1-q^{\alpha}}} & {\rm if}\; \gamma =0,\; q > 1;
\\ \smallskip \\
\lambda _n \ra \textstyle{\frac {q^{\beta}}{1-q^{\gamma}}}
                       & {\rm if}\; \gamma >0,\; \alpha=0,\; 0< q <1;
\ea\eeq
and if $n\to -\infty$ we have
\begin{equation}\label{f20a1a}\begin{array}{ll}
\lambda _n \ra -\infty & {\rm if}\;\cases{
		    \gamma =0,\;    q > 1,  & or \cr
		    \gamma >0,\; 0< q < 1,  & or \cr
         \gamma >0,\;\alpha <0,\;    q > 1;  &    \cr}
\\ \smallskip \\
\lambda _n \ra 0^-     & {\rm if}\; \gamma >0,\;\alpha >0,\;  q > 1;
\\ \smallskip \\
\lambda _n \ra \lambda _0- \textstyle{\frac {q^{\alpha \nu _0 +\beta}}
       {q^{\alpha} -1}} & {\rm if}\; \gamma =0,\; 0 < q < 1;
\\ \smallskip \\
\lambda _n \ra \textstyle{\frac {q^{\beta}}{1-q^{\gamma}}}
                       & {\rm if}\; \gamma >0,\; \alpha=0,F\;  q > 1.
\ea\eeq

Denote by $n_0$ the greatest integer for which in (\ref{f16}) we have
$\lambda _n\leq 0.$ In this case we inevitably have $\psi _{n_0}=0$
(because of $a^{+}a\geq 0).$ Then $$a\psi _{n_0+1}=\psi _{n_0}=0,\quad
{\hbox{\rm where}}\quad \psi _{n_0+1}\neq 0,$$
and the vector $\psi _{n_0+1}$ is the common eigenvector for the operators
$a^{+}a$ and $N$ with eigenvalues $\lambda _{n_0+1}=0$ and
$\nu _0^{\prime }=\nu _0+n_0+1,$ respectively. Moreover the following
relation $$\lambda _0+q^{\alpha \nu _0+\beta }\frac{q^{(n_0+1)
(\alpha -\gamma )}-1} {q^\alpha -q^\gamma }=0.$$ is hold.

Now we repeat the whole construction above with the following substitutions
$$
\psi _0\rightarrow \psi _{n_0+1},\qquad \lambda _0\rightarrow 0,\qquad \nu
_0\rightarrow \nu _0^{\prime }.
$$
Using for the normalized basis states notation $|n>:=\psi
_{n_0+1+n},\;n=0,1,\ldots $ ($n\geq 0$) we receive the representation $\pi
(\nu _0^{\prime }|{\cal W}_{\alpha ,\beta }^\gamma )$ of the algebra
${\cal W}_{\alpha ,\beta }^\gamma (q)$ in the form

\begin{equation}
\label{f21}\left\{
\begin{array}{lcc}
a^{+}|n>=q^{\alpha \nu _0^{\prime }/2}\left( F_{\alpha ,\beta }^\gamma
(n+1)\right) ^{1/2}|n+1>, & \quad  &  \\
a|n>=q^{\alpha \nu _0^{\prime }/2}\left( F_{\alpha ,\beta }^\gamma
(n)\right) ^{1/2}|n-1>, &  & n\geq 0, \\
N|n>=(\nu _0^{\prime }+n)|n>. &  &
\ea \right.
\end{equation}
These series of irreducible representations of quasi-Fock type are labeled by
the real number $\nu _0^{\prime }$ which gives the lowest bound of the
spectrum of the operator $N.$ For $\nu _0^{\prime }=0$ we receive the q-Fock
representation. {\bf Only }in this q-Fock representation together with (\ref
{f1}) additional relation
\begin{equation}
\label{f22}aa^{+}-q^\alpha a^{+}a=q^{\gamma N+\beta }
\end{equation}
is hold.

According to (\ref{f22}) for the representations
$\pi (\nu _0^{\prime }|{\cal W}_{\alpha ,\beta }^\gamma )$ (\ref{f21})
the element $K$ (\ref{f8}) fulfils the inequality
\begin{equation}
\label{f23}K=F_{\alpha ,\beta }^\gamma (N+I)-q^\alpha F_{\alpha ,\beta
}^\gamma (N)=q^{\gamma N+\beta }>0, \end{equation}
as it indeed must be because the representation
$\pi (\nu _0^{\prime }|{\cal W}_{\alpha ,\beta }^\gamma )$ is non-singular
if we extend for the case of
the algebra ${\cal W}_{\alpha ,\beta }^\gamma (q)$ criterion developed in
\cite{34} for more familiar algebra ${\cal W}_{-1,0}^1(q).$

We also note in closing of the consideration of this case that for
$0<q<1,\;0<\alpha <\gamma$ operators $a^+a,\,aa^+ \in${\got S}$_1 $ but if
$0<q<1,\;0=\alpha <\gamma$ or if $q>1,\;\alpha <\gamma =0$ these operators
are bounded.

In the next possible case
$${\bf (B)\qquad }\quad 0<q<1,\;\quad \lambda _0>\frac{q^{\alpha
(\nu _0-1)+\beta }}{1-q^{\gamma -\alpha }}$$
 we consider first the case $ \alpha < 0.$ Then
according (\ref{f16}) $\lambda _n\rightarrow \infty $ as $|n|\ra \infty .$
Therefore it follows that the operator $a^{+}a$ has the lowest
eigenvalue for some value of $n_0.$ Now
we can repeat all of the considerations above starting from the values
$\lambda _{n_0},\,\psi _{n_0}$ and $\nu _0+n_0$ in place of $\lambda
_0,\,\psi _0$ and $\nu _0,$ respectively. Then we once more obtain
the relation (\ref{f16}) with an additional condition $\lambda _n\geq
\lambda _0.$

For $n>0$ this inequality gives
$$\lambda _0<q^{\alpha \nu _0}\frac{F_{\alpha ,\beta }^\gamma (n)}
{1-q^{n\gamma }},$$ and for $n<0$ we obtains
$$\lambda _0>-q^{\alpha \nu _0}\frac{F_{\alpha ,\beta }^\gamma (n)}
{q^{n\gamma }-1}.$$ This both inequalities are hold (with new values
of $\lambda _0,\nu _0)$ if \begin{equation}\label{f24}-q^{\alpha \nu _0}
\frac{F_{\alpha ,\beta }^\gamma (-1)}{q^{-\gamma }-1}<\lambda _0<q^{\alpha
\nu _0}\frac{F_{\alpha ,\beta }^\gamma (1)}{1-q^\gamma }.\end{equation}

As in the reference \cite{7}, in this case we have two-parameter family of
non equivalent irreducible representations $\pi (\nu _0,\lambda _0|
{\cal W}_{\alpha ,\beta }^\gamma ).$ On the elements of the orthogonal
basis $$ \left\{ |n>= \psi _{n_0+n}|n\in {\hbox{\BBm Z}}\right\} $$
generators of the algebra act according the formulae
\begin{equation}
\label{f25}\left\{
\begin{array}{lcc}
a^{+}|n>=\left( q^{(n+1)\gamma }\lambda _0+q^{\alpha \nu _0}F_{\alpha ,\beta
}^\gamma (n+1)\right) ^{1/2}|n+1>, & \quad  &  \\
a\;|n>=\left( q^{n\gamma }\lambda _0+q^{\alpha \nu _0}F_{\alpha ,\beta
}^\gamma (n)\right) ^{1/2}|n-1>, &  & n\in {\hbox{\BBm Z}}, \\
N\;|n>=(\nu _0+n)|n>. &  & \end{array}\right.\end{equation}
The irreducible representations of this family $\pi (\nu _0,\lambda _0|
{\cal W}_{\alpha ,\beta }^\gamma )$ are singular ones, they disappears
in the ''classical limit'' $q\rightarrow 1$ . For these representations
\begin{equation} \label{f26} \begin{array}{l}
K=\left( q^{(N+I)\gamma }\lambda _0+q^{\alpha \nu _0}F_{\alpha ,\beta
}^\gamma (N+I)\right) -q^\alpha \left( q^{N\gamma }\lambda _0+q^{\alpha \nu
_0}F_{\alpha ,\beta }^\gamma (N)\right) = \\
\qquad \;\,q^{N\gamma }\lambda _0(q^\gamma -q^\alpha )+q^{\alpha \nu
_0}q^{\gamma N+\beta }=\,q^{N\gamma }(q^\gamma -q^\alpha )\left(
\lambda _0 - \frac{q^{\alpha \nu _0+\beta }}{q^\gamma -q^\alpha }\right) <0.
\end{array}
\end{equation}
[In deriving the relation (\ref{f26}) we take into account conditions
$\lambda _0>\frac{q^{\alpha (\nu _0-1)+\beta }}{1-q^{\gamma -\alpha }},$ and
$0<q<1,$ which means that $\alpha <\gamma \Rightarrow q^\alpha >q^\gamma ].$

We note that the singularity of the representation is also indicated by the
fact that the inequalities  (\ref{f24}) became contradictory in the
''classical limit'' $q\rightarrow 1$ .

In the opposite case $\alpha \leq 0\ $ according to  (\ref{f16}) we have
$\lambda _n\rightarrow +\infty $ as $n\rightarrow -\infty $ and
$$
\lambda_n\rightarrow 0^{+},\qquad {\rm if}\quad\alpha<0;\quad\
\lambda_n\rightarrow \frac{q\beta}{1-q^\gamma}\qquad {\rm if}\quad \alpha=0,
$$
as $n\rightarrow +\infty$. Then we receive again the representations
$\pi (\nu _0,\lambda _0|
{\cal W}_{\alpha ,\beta }^\gamma )$(\ref{f25})
without the restriction (\ref{f24})
The inequality (\ref{f26}) is hold so these representations are singular ones.
Now we must consider the last possible case
$$
{\bf (C)\qquad }\quad 0<q<1,\;\quad \lambda _0=\frac{q^{\alpha (\nu
_0-1)+\beta }}{1-q^{\gamma -\alpha }}.
$$
In this case condition (\ref{f16}) gives

\begin{equation}
\label{f27}\quad \lambda _n=\frac{q^{\alpha (\nu _0+n)+\beta }}{q^\alpha
-q^\gamma },\qquad n\in {\hbox{\BBm Z}},\end{equation} which means that
there exists no minimal as well as maximal eigenvalue for $a^{+}a.$ In the
irreducible representations
$\pi _S(\nu _0|{\cal W}_{\alpha ,\beta }^\gamma )$
of this one-parameter family the action of generators is given by
\begin{equation}
\label{f28}\left\{
\begin{array}{lcc}
a^{+}|n>=q^{\frac{\alpha \nu _0+\beta +\alpha (n+1)}2}\left( q^\alpha
-q^\gamma \right) ^{-1/2}|n+1>, & \quad  &  \\
a\;|n>=q^{\frac{\alpha \nu _0+\beta +\alpha n}2}\left( q^\alpha -q^\gamma
\right) ^{-1/2}|n-1>, &  & n\in {\hbox{\BBm Z}}, \\
N\;|n>=(\nu _0+n)|n>. &  &
\end{array}
\right.
\end{equation}
In the representations $\pi _S(\nu _0|{\cal W}_{\alpha ,\beta }^\gamma )$ of
this one-parameter family additional relation
\begin{equation}
\label{f29}aa^{+}-q^\alpha a^{+}a=0
\end{equation}
is hold. This condition means that in such representations $K=0.$The
representations $\pi _S(\nu _0|{\cal W}_{\alpha ,\beta }^\gamma )$ are
singular, and called strange [special] ones because in ''classical limit''
$q\rightarrow 1$ commutation relations (\ref{f29}) and (\ref{f1}) became
contradictory.

Now we briefly consider the second case from the list (\ref{f5}).

2) The case $\gamma \geq 0,\quad \alpha >\gamma $.

As above (see (\ref{f20a}) ) we suppose that the inequalities
$$
{\bf (A)\qquad }0<q<1\quad {\rm or}\quad q>1,\;\lambda _0<\frac{q^{\alpha
(\nu _0-1)+\beta }}{1-q^{\gamma -\alpha }}.
$$
are hold and by the same reasoning as in the case 1) we received the
nonsingular representation (\ref{f21}).  In this case we also have $K>0.$
Note that for this case we have
$$
\begin{array}{lll}
\begin{array}{l}
\lambda _n\rightarrow 0^{+}, \\
\lambda _n\rightarrow -\infty ,
\end{array}
&
\begin{array}{l}
{\rm if\;n\rightarrow +\infty } \\ {\rm if\;n\rightarrow -\infty }
\end{array}
& {\rm for}\;0<q<1; \\
\begin{array}{l}
\lambda _n\rightarrow +\infty , \\
\lambda _n\rightarrow 0^{-},
\end{array}
&
\begin{array}{l}
{\rm if\;n\rightarrow +\infty } \\ {\rm if\;n\rightarrow -\infty }
\end{array}
& {\rm for}\;q>1.
\end{array}
$$
It can be easely checked that in the case $0<q<1,$ operators $a^{+}a$
and $aa^{+}$ are compact (or from the class {\got S}$_{\infty}$ ) and
even nuclear ones (or from the class {\got S}$_{1}$ ).

In the next case
$$
{\bf (B)\qquad }\quad q>1,\;\quad \lambda _0>\frac{q^{\alpha (\nu
_0-1)+\beta }}{1-q^{\gamma -\alpha }}.
$$
we have the representation (\ref{f25}). The action of operator $K$ on
vector $|n>$ gives
\begin{equation}
\label{f30}K|n>=\left( aa^{+}-q^\alpha a^{+}a\right) |n>=q^{\gamma
n}(q^\alpha -q^\gamma )\left( \frac{q^{\alpha \nu _0+\beta }}{q^\alpha
-q^\gamma }-\lambda _0\right) |n>
\end{equation}
so that $K<0,$ because in this case
$(q^\alpha -q^\gamma )>0$
and $\left( \frac{q^{\alpha \nu _0+\beta }}{q^\alpha -q^\gamma }-\lambda
_0\right) <0$. Thus this representations are singular ones. We have here
$$
\begin{array}{lll}
\begin{array}{l}
\lambda _n\rightarrow +\infty , \\
\lambda _n\rightarrow 0^{+},
\end{array}
&
\begin{array}{l}
{\rm if\;n\rightarrow +\infty, } \\ {\rm if\;n\rightarrow -\infty .}
\end{array}
& {\rm (}q>1),
\end{array}
$$
If we take into account the condition on $\lambda _0$ and note that
$\lambda_n\geq 0$ and $F_{\alpha ,\beta }^\gamma (n)<0$  if $n<0$ then
we have
$$
0\leq \lambda _n\leq \lambda _0\quad {\rm if\quad }n<0
$$
for such representations.

In the last possible case
$$
{\bf (C)\qquad }\quad q>1,\;\quad \lambda _0=\frac{q^{\alpha (\nu
_0-1)+\beta }}{1-q^{\gamma -\alpha }},
$$
for all values $n\in {\hbox{\BBm Z}}$ we have
\begin{equation}\label{f31}\lambda _n=
\frac{q^{\alpha (\nu _0+n)+\beta }}{q^\alpha -q^\gamma }\geq 0,
\end{equation}
and

$$
\begin{array}{lll}
\begin{array}{l}
\lambda _n\rightarrow +\infty , \\
\lambda _n\rightarrow 0^{+},
\end{array}
&
\begin{array}{l}
{\rm if\;n\rightarrow +\infty, } \\ {\rm if\;n\rightarrow -\infty. }
\end{array}
& {\rm (}q>1).
\end{array}
$$
As in the case 1) we received the representations of special type
(\ref{f28}) for which $K=0,$ that is representations are  singular and
additional relation (\ref{f29}) is hold too. In representations of
this type the operator $a^{+}a$ has no minimal eigenvalue.

3) The case $\gamma \leq 0,\quad \alpha <\gamma $.

$$
{\bf (A)\qquad }q>1\quad {\rm or}\quad 0<q<1,\;\lambda _0<\frac{q^{\alpha
(\nu _0-1)+\beta }}{1-q^{\gamma -\alpha }}.
$$

In this subcase $K>0$ and we have nonsingular representations (\ref{f21})
for which

$$
\begin{array}{lll}
\begin{array}{l}
\lambda _n\rightarrow 0^{+}, \\
\lambda _n\rightarrow -\infty ,
\end{array}
&
\begin{array}{l}
{\rm if\;n\rightarrow +\infty } \\ {\rm if\;n\rightarrow -\infty }
\end{array}
& {\rm for}\;q>1; \\
\begin{array}{l}
\lambda _n\rightarrow +\infty , \\
\lambda _n\rightarrow 0^{-},
\end{array}
&
\begin{array}{l}
{\rm if\;n\rightarrow +\infty } \\ {\rm if\;n\rightarrow -\infty }
\end{array}
& {\rm for}\;0<q<1.
\end{array}
$$
As in the case 2A) (but now for $q>1$) operators $a^{+}a$  and
$aa^{+}$ are nuclear (from {\got S}$_1$-class).

$$
{\bf (B)\qquad }\quad 0<q<1,\;\quad \lambda _0>\frac{q^{\alpha (\nu
_0-1)+\beta }}{1-q^{\gamma -\alpha }}
$$
In this subcase representations have the form (\ref{f25}), $K<0$ and
representations are singular.

$$
{\bf (C)\qquad }\quad 0<q<1,\;\quad \lambda _0=\frac{q^{\alpha (\nu
_0-1)+\beta }}{1-q^{\gamma -\alpha }}.
$$
The inequality (\ref{f31}) holds for all values $n\in {\hbox{\BBm Z}}.$
In this case
$K=0,$ and the representations are singular and  have the form
(\ref{f28}). The additional relation   (\ref{f29}) also hold in this case.
We have not  the minimal eigenvalue for $a^{+}a$ and
$$
\begin{array}{lll}
\begin{array}{l}
\lambda _n\rightarrow +\infty , \\
\lambda _n\rightarrow 0^{+},
\end{array}
&
\begin{array}{l}
{\rm if\;n\rightarrow +\infty ,} \\ {\rm if\;n\rightarrow -\infty. }
\end{array}
& {\rm (0<}q<1).
\end{array}
$$

4) The case $\gamma \leq 0,\quad \alpha >\gamma $.

$$
{\bf (A)\qquad }0<q<1\quad {\rm or}\quad q>1,\;\lambda _0<\frac{q^{\alpha
(\nu _0-1)+\beta }}{1-q^{\gamma -\alpha }}.
$$
In this subcase $K>0$ and we have nonsingular representations (\ref{f21}).
Note that now we have

$$
\begin{array}{lll}
\begin{array}{l}
\lambda _n\rightarrow +\infty , \\
\lambda _n\rightarrow 0^{+}, \\
\lambda _n\rightarrow q^{\alpha \nu _0+\beta }
\end{array}
&
\begin{array}{l}
{\rm if\;n\rightarrow +\infty } \\ {\rm if\;n\rightarrow +\infty } \\
{\rm  if\;n\rightarrow +\infty }
\end{array}
&
\begin{array}{l}
{\rm for}\;0<q<1,\;{\rm or}\;q>1\;{\rm and\;\alpha }>0; \\ {\rm for}\;q>1\;
{\rm and\;\alpha }<0; \\ {\rm for}\;q>1\;{\rm and\;\alpha }=0;
\end{array}
\\
\begin{array}{l}
\lambda _n\rightarrow 0^{-}, \\
\lambda _n\rightarrow -\infty , \\
\lambda _n\rightarrow -q^{\alpha \nu _0+\beta }
\end{array}
&
\begin{array}{l}
{\rm if\;n\rightarrow -\infty } \\ {\rm if\;n\rightarrow -\infty } \\
{\rm if\;n\rightarrow -\infty }
\end{array}
&
\begin{array}{l}
{\rm for}\;0<q<1\;{\rm and\;\alpha }<0; \\ {\rm for}\;0<q<1\;{\rm
and\;\alpha }>0,\;{\rm or}\;q>1; \\ {\rm for}\;0<q<1\;
{\rm and\;\alpha =}0;
\end{array}
\end{array}
$$
This means that if $q>1\;{\rm and\;\alpha }<0,$ then operators $a^{+}a$
and $aa^{+}$ are nuclear (from {\got S}$_1$-class), and
${\rm for}\;q>1\;{\rm and\;\alpha }=0$
this operators are bounded.

The cases

$$
{\bf (B)\qquad }\quad q>1,\;\quad \lambda _0>\frac{q^{\alpha (\nu
_0-1)+\beta }}{1-q^{\gamma -\alpha }}
$$
and

$$
{\bf (C)\qquad }\quad q>1,\;\quad \lambda _0=\frac{q^{\alpha (\nu
_0-1)+\beta }}{1-q^{\gamma -\alpha }}.
$$
are completely analogous to the cases 1B) and 1C). In the last case

5) $\qquad \alpha =\gamma $

for all values $0<q<\infty $ we have only nonsingular  representations of
the type (\ref{f21}) with $K>0.$ In these representations

$$
\begin{array}{lll}
\begin{array}{l}
\lambda _n\rightarrow 0^{+}, \\
\lambda _n\rightarrow +\infty ,
\end{array}
&
\begin{array}{l}
{\rm if\;n\rightarrow +\infty } \\ {\rm if\;n\rightarrow +\infty }
\end{array}
&
\begin{array}{l}
{\rm for}\;0<q<1,\,\gamma >0\;{\rm or\;}q>1,\,\gamma <0; \\ {\rm for}
\;q>1,\,\gamma >0\;{\rm or\;}0<q<1,\,\gamma \leq 0;
\end{array}
\\
\begin{array}{l}
\lambda _n\rightarrow -\infty , \\
\lambda _n\rightarrow 0^{-},
\end{array}
&
\begin{array}{l}
{\rm if\;n\rightarrow -\infty } \\ {\rm if\;n\rightarrow -\infty }
\end{array}
&
\begin{array}{l}
{\rm for}\;0<q<1,\,\gamma >0\;{\rm or\;}q>1,\,\gamma \leq 0; \\ {\rm for}
\;q>1,\,\gamma >0\;{\rm or\;}0<q<1,\,\gamma <0.
\end{array}
\end{array}
$$
Thus in the cases $\;0<q<1,\,\gamma >0\;{\rm or\;}q>1,\,\gamma <0$ we have
nuclear ({\got S}$_1$-class) operators $a^{+}a$ and $aa^{+}.$

{\bf Remarks}. 1. Let us consider the modification of the spectral
properties of the operators $a^{+}a$ ( and $aa^{+}$) caused by changing
of the value of the deformation parameter $q$ near the classical point
$q=1$ for some concrete variants of the algebras ${\cal W}^{\gamma}_{\alpha,
\beta}(q)$ mentioned in the section 2. For the Tamm-Dancoff algebra
${\cal W}^{1}_{1,0}(q)$ operators $a^{+}a,\; aa^{+}$ are nuclear (that is of
the {\got S}$_1$ class) if $0<q<1$ (see case 5) in the section 3). But when
$q\geq 1$ the eigenvalues $\lambda _n$ of the operator $a^{+}a$ tending to
infinity ($\lambda _n\ra +\infty$) as $n\ra +\infty$. The same result is also
true for the operator $aa^{+}$ ($\mu _n \ra +\infty$ as $n\ra +\infty$).

In the case of the Arik-Coon-Kuryskin algebra ${\cal W}^{1}_{0,0}(q)$ for
the case $0<q<1$ and $\lambda _0 <{\frac 1{1-q}}$ (see case A) in the
subsection 3.01) operators $a^{+}a$ and $aa^{+}$ are bounded with the norm
growing as $(1-q)^{-1}$ as $q\ra 1^-$. When $q\geq 1$ than $\lambda _n\ra +
\infty$ and $\mu _n \ra +\infty$  if $n\ra +\infty$.

For the algebra ${\cal W}^{1}_{-1,0}(q)$ there are no quality changes in
spectral properties of the operators $a^{+}a$ ( and $aa^{+}$) in the process
of passing the point $q=1$.

For the Feinsilver case ${\cal W}^{0}_{-2,0}(q)$ considered operators are
bounded when $q>1$ with the norm growing as $(q-1)^{-1}$ as $q\ra 1^+$ and
if $0<q<1$ and  $\lambda _0 <{\frac {q^{-2(\nu _0 -1)}}{1-q}}$ then
$\lambda _n\ra +\infty$ and $\mu _n \ra +\infty$ as $n\ra +\infty$.

The mentioned peculiarities in spectral properties may be helpful in
attempts of construction of physical models with help of deformed oscillator
operators.

2.  As mentioned above some authors together with the main relations
defining the algebra ${\cal W}^{\gamma}_{\alpha,\beta}(q)$
\beq\label{f32}aa^+-q^{\gamma}a^+a=q^{\alpha N+\beta} \eeq
required also additional ones, which in general have the form
\beq\label{f33}aa^+-q^{\tilde\gamma}a^+a=q^{{\tilde{\alpha}}N+
{\tilde{\beta}}}.\eeq
Because from the equality $\gamma = \tilde{\gamma}$ it follows that
$\alpha=\tilde{\alpha}$ we can restrict our attention to the case
$\gamma \neq \tilde{\gamma}$.

If both of the above relations (\ref{f32}) and (\ref{f33}) are hold then we
have
\beq\label{f34}\ba
a^{+}a=(q^{-\gamma}-q^{-\tilde\gamma})^{-1}(-q^{\tilde\alpha N +\tilde\beta
-\gamma-\tilde\gamma}+q^{\alpha N +\beta -\gamma-\tilde\gamma}) \\
aa^{+}=(q^{-\gamma}-q^{-\tilde\gamma})^{-1}(q^{\alpha N +\beta -\gamma}
-q^{\tilde\alpha N +\tilde\beta -\tilde\gamma}) \ea \eeq
and
\beq\label{f35}
K=aa^{+}-q^{\alpha}a^{+}a=(q^{-\gamma}-q^{-\tilde\gamma})^{-1}
\bigl( q^{\alpha N +\beta -\gamma}(1-q^{\alpha -\tilde\gamma})-
q^{\tilde\alpha N +\tilde\beta -\tilde\gamma}(1-q^{\alpha -\gamma})\bigr).
\eeq
Note that only nonsingular representations appears in the two most popular
cases:

1) The $q\leftrightarrow q^{-1}$-symmetrical case under which (\ref{f33}) can
be obtained from (\ref{f32}) after such substitution (thus
$\alpha =- \tilde\alpha$, $\gamma =- \tilde\gamma$ and $\beta =-
\tilde\beta$).

2) The case when (\ref{f33}) can be obtained from (\ref{f32}) after the
substitution $\alpha \leftrightarrow \gamma $ ($\alpha =\tilde\gamma$
and $\tilde\alpha =\gamma$).

\section{The Hopf algebra structure}

The Hopf algebra structure for the quantum deformed oscillator algebras
${\cal W}^{\gamma}_{\alpha,\beta}(q)$ is not known, and probably non exist.
But for some modificated algebras such structure can be introduced
as it was shown in \cite{63}. In this section we
consider two examples of such algebras admiting the Hopf algebra
structure using the similar argumentation.

(1) From the relation
\beq\label{f36}aa^+-q^{\gamma}a^+a=q^{\alpha N+\beta} \qquad \alpha\neq 0,\;
\gamma \neq 0, \eeq
and its $q\leftrightarrow q^{-1}$ partner
\beq\label{f37}aa^+-q^{-\gamma}a^+a=q^{-\alpha N-\beta}.  \eeq
we received
\beq\label{f38}\ba
aa^+= F^{1}_{-1,0}(-{\frac {\alpha}{\gamma}}N-{\frac {\beta}{\gamma}}+1;
q^{-\gamma}) \\  \smallbreak \\
a^+a=F^{1}_{-1,0}(-{\frac {\alpha}{\gamma}}N-{\frac {\beta}{\gamma}};
q^{-\gamma}). \ea \eeq
and
\beq\label{f39}[a,a^+]=
F^{1}_{-1,0}(-{\frac {\alpha}{\gamma}}N-
{\frac {\beta}{\gamma}}+1;q^{-\gamma})-
F^{1}_{-1,0}(-{\frac {\alpha}{\gamma}}N-{\frac {\beta}{\gamma}};q^{-\gamma}).
\eeq

Let us consider the oscillator algebra
${\cal W}^{\gamma}_{\alpha,\beta}(q;q^{-1})$, defined by this relation and
equations of motion (\ref{f2}). Then the coproduct $\Delta$, counity
$\varepsilon$ and antipod $S$ in ${\cal W}^{\gamma}_{\alpha,\beta}(q;q^{-1})$
can be defined by the relations
\begin{equation}\label{f40}\begin{array}{c}\begin{array}{c}
\Delta (a^{+})=\!c_1a^{+}\otimes q^{\alpha _1\,N}+c_2q^{\alpha _2N}\otimes
a^{+};\quad \Delta (a)=c_3a\otimes q^{\alpha _3\,N}+c_4q^{\alpha _4N}\otimes
a; \\ \Delta (N)=\!c_5N\otimes {\bf 1}+c_6{\bf 1}\otimes N+\gamma _1{\bf
1\otimes 1}; {\bf \quad }\Delta ({\bf 1})={\bf 1\otimes 1}; \end{array} \\
\varepsilon (a^{+})=\!c_7;\quad \varepsilon (a)=c_8;\quad \varepsilon
(N)=c_9;\quad \varepsilon ({\bf 1})=1; \\ S(a^{+})=\!-c_{10}a^{+};\quad
S(a)=-c_{11}a;\quad S(N)=-c_{12}N+c_{13}{\bf 1};\quad S({\bf 1})=1;
\end{array} \end{equation}
The constants $c_i,\alpha _k$ and $\gamma _1$can be easely finded from
concordance conditions
\begin{equation}\label{f41}\begin{array}{c} ({\rm id}\otimes \Delta )\Delta
(h)=(\Delta \otimes {\rm id})\Delta (h), \\ ({\rm id}\otimes \varepsilon )
\Delta (h)=(\varepsilon \otimes {\rm id})\Delta (h), \\ m({\rm id}
\otimes S)\Delta (h)=m(S\otimes {\rm id})\Delta (h)=\varepsilon (h){\bf 1},
\end{array} \end{equation}
for all $h\in {\cal H,}$ where ${\cal H}$ is a considered Hopf algebra and
$m:{\cal H} \rightarrow {\cal H\otimes H}$ is the multiplication in ${\cal
H}$. From (\ref{f40}), (\ref{f41}) and (\ref{f2}) one obtains
\begin{equation} \label{f42} \begin{array}{c}
c_1=q^{\alpha _1\,\gamma _1};\quad c_2=q^{\alpha _2\,\gamma _1};\quad
c_3=q^{\alpha _3\,\gamma _1};\quad c_4=q^{\alpha _4\,\gamma _1};\quad \\
c_5=c_6=1;\quad c_7=c_8=0;\quad c_9=-\gamma _1;\quad c_{10}=q^{\alpha
_1\,};\quad \\ c_{11}=q^{-\alpha _3};\quad c_{12}=1;\quad c_{13}=-2
\gamma _1;\quad \alpha _2=-\alpha _1;\quad \alpha _4=-\alpha _3.\quad
\end{array} \end{equation}
Moreover $\Delta $ and $\varepsilon $ are homomorphisms and $S$ is
antihomomorphism of the deformed algebra. In particular we must have
\begin{equation}
\label{f43}\Delta (a)\Delta (a^{+})-\Delta (a^{+})\Delta (a)\!=\Delta
(F^{1}_{-1,0}(-{\frac {\alpha}{\gamma}}N-
{\frac {\beta}{\gamma}}+1;q^{-\gamma})\!-\!
F^{1}_{-1,0}(-{\frac {\alpha}{\gamma}}N-{\frac {\beta}{\gamma}};q^{-\gamma}).
\end{equation}
This gives the relations
\begin{equation}
\label{f44}q^{\alpha _1\,}=q^{\alpha _3};\quad \quad \alpha _1+\alpha
_3=\alpha ;\quad \quad q^{2\alpha \,\gamma _1}=-q^{2\beta -\gamma }.
\end{equation}
For real $q$ (\ref{f44}) imply that
\begin{equation}
\label{f45}\alpha _1=\alpha _3=\frac \alpha 2;\quad \quad
\gamma _1=\frac{
2\beta -\gamma }{2\alpha }+i\frac{(2k+1)\pi }{2\alpha \ln q},\;k\in
{\hbox{\BBm Z}},\end{equation}
which fix all constants in (\ref{f40}). Thus we received finally the
following definition of Hopf algebra structure in concidered case
\begin{equation}
\label{f46}
\begin{array}{c}
\begin{array}{c}
\Delta (a^{+})=a^{+}\otimes q^{\alpha (N+\gamma _1)/2}+q^{-\alpha (N+\gamma
_1)/2}\otimes a^{+}; \\
\Delta (a)=a\otimes q^{\alpha (N+\gamma _1)/2}+q^{-\alpha (N+\gamma
_1)/2}\otimes a; \\
\Delta (N)=N\otimes {\bf 1}+{\bf 1}\otimes N+\gamma _1{\bf 1\otimes 1};{\bf
\quad }\Delta ({\bf 1})={\bf 1\otimes 1};
\end{array}
\\
\varepsilon (a^{+})=0=\varepsilon (a);\quad \varepsilon (N)=-\gamma _1;\quad
\varepsilon (
{\bf 1})=1; \\ S(a^{+})=-q^{-\alpha /2}a^{+};\quad S(a)=-q^{-\alpha
/2}a;\quad S(N)=-N-2\gamma _1{\bf 1};\quad S({\bf 1})=1.
\end{array}
\end{equation}

2). Consider the deformed oscillator algebra generated by the relations
(\ref{f2}) and
\begin{equation}\label{f47}[a,a^{+}]=F_{-\gamma ,\beta }^\gamma (N+1;q)-
F_{-\gamma ,\beta}^\gamma (N;q).\end{equation}
Repeating the calculation we obtain relations (\ref{f46}) but with $\alpha $
replaced by $\gamma $ and $
\gamma _1={\frac 12}-i\frac{(2k+1)\pi }{2\gamma \ln q},\;k\in${\BBm Z}. We
recall finally that relations (\ref{f47}) follows from couple of commutation
rules \begin{equation} \label{ff1}aa^{+}-q^\gamma a^{+}a=q^{-\gamma N+\beta
},\quad aa^{+}-q^{-\gamma }a^{+}a=q^{\gamma N+\beta }.  \end{equation}

\noindent{\bf Acknowledgments.}
\noindent
The authors would like to thank P.P.Kulish for helphul comments and
valuable discussions of numberous subjects related with the theme of this
work.

V.B. thanks A. Laptev, T. Weidl and the Royal Institute of Technology
Stockholm for hospitality and support.

The work of E.D. and S.E. was supported by the Russian
Foundation for Fundamental Research under the Grant N 95-01-00569-a.
\vspace{.5cm}


\end{document}